\documentclass[12pt]{amsart}
\usepackage{amsmath}
\usepackage{amscd}
\usepackage{amssymb}
\usepackage{float}

 \usepackage[all,cmtip]{xy}
\usepackage{graphicx}
 \usepackage{placeins}
 \usepackage{float}
\usepackage{color}
\usepackage{amstext}    
\usepackage{array}

\newtheorem{thm}{Theorem}[section]

\newtheorem{lemma}[thm]{Lemma}

\newtheorem{proposition}[thm]{Proposition}

\newtheorem{thm-dfn}[thm]{Theorem-Definition}

\theoremstyle{definition}

\newtheorem{remark}[thm]{Remark}

\numberwithin{equation}{section}

\newcommand{\cB}{{\mathcal B}}

\newcommand{\cO}{{\mathcal O}}
\newcommand{\cA}{{\mathcal A}}

\newcommand{\cN}{{\mathcal N}}

\newcommand{\cS}{{\mathcal S}}
\newcommand{\cE}{{\mathcal E}}
\newcommand{\cU}{{\mathcal U}}

\newcommand{\bC}{{\mathbb C}}
\newcommand{\bZ}{{\mathbb Z}}

\newcommand{\bQ}{{\mathbb Q}}

\newcommand{\Lb}{{\mathfrak{b}}}
\newcommand{\Lt}{{\mathfrak{t}}}
\newcommand{\Lg}{{\mathfrak g}}

\newcommand{\Ln}{{\mathfrak{n}}}

\newcommand{\Ll}{{\mathfrak{l}}}

\newcommand{\tk}{{\textbf{k}}}

\newcommand{\on}{\operatorname}

\newcommand{\beqn}{\begin{equation*}}
\newcommand{\eeqn}{\end{equation*}}
\newcommand{\beq}{\begin{equation}}
\newcommand{\eeq}{\end{equation}}

\newcommand{\N}{\cN_{\Lg^*}}

\begin{document}

\title[Springer correspondence for  exceptional Lie algebras]
{Springer correspondence for exceptional Lie algebras and their duals in small characteristic}
        \author{Ting Xue}
        \address{School of Mathematics and Statistics, University of Melbourne, Australia and Department of Mathematics and Statistics, University of Helsinki, Finland}
        \email{ting.xue@unimelb.edu.au}
        \thanks{The author was supported in part by the Academy of Finland and the ARC grant DE160100975 and  DP150103525.}

\maketitle
\begin{abstract}
We describe the Springer correspondence explicitly for exceptional Lie algebras of type $G_2$ and $F_4$ and their duals in bad characteristics, i.e. in characteristics 2 and 3. 
\end{abstract}
\section{Introduction}
In this paper we describe the Springer correspondence explicitly for exceptional Lie algebras of type $G_2$ and $F_4$ and their duals in bad characteristics. In what follows we describe the set-up of the paper in more detail. 

Let $G$ be an almost-simple algebraic group defined over an algebraically closed field $\tk$ of prime characteristic and $\Lg$ the Lie algebra of $G$. Let $\Lg^*$ be the dual vector space of $\Lg$.  Denote by $\cU_G$, $\cN_\Lg$ and $\cN_{\Lg^*}$ the varieties of unipotent elements in $G$, nilpotent elements in $\Lg$, and nilpotent elements in $\Lg^*$, respectively. We recall that an element in $\Lg^*$ is called nilpotent if it annihilates a Borel subalgebra of $\Lg$ (see \cite{KW}). The group $G$ acts on $\Lg$ and $\Lg^*$ by adjoint and coadjoint action, respectively. It is known that the number of $G$-orbits in $\cU_G$, $\cN_\Lg$ and $\cN_{\Lg^*}$, respectively, is finite (for $\cN_{\Lg^*}$ see \cite{X1} and references there). We fix a prime $l\neq\on{char}\tk$. Let $\cA_G$ (resp. $\cA_{\Lg}$, $\cA_{\Lg^*}$) denote the set of all pairs $(\cO,\cE)$, where $\cO\subset\cU_G$ (resp. $\cO\subset\cN_\Lg$, $\cO\subset\N$) is a $G$-orbit and $\cE$ is an irreducible $G$-equivariant $\bar\bQ_l$-local system on $\cO$ (up to isomorphism). Let $W$ be the Weyl group of $G$ and let $\on{Irr}(W)$ denote the set of irreducible characters (over $\bar\bQ_l$) of the Weyl group $W$. The Springer correspondence for $G$ (resp. $\Lg$, $\Lg^*$) maps the set $\on{Irr}(W)$ injectively into the set $\cA_G$ (resp. $\cA_{\Lg}$, $\cA_{\Lg^*}$); we denote the map by $\gamma_G$ (resp. $\gamma_\Lg$, $\gamma_{\Lg^*}$). The construction was originally given for $\Lg$ by Springer in large characteristic \cite{Sp} and generalized to arbitrary characteristic for $G$ by Lusztig \cite{L3}. The same construction as Lusztig's works for $\Lg$ (resp. $\Lg^*$) in any characteristic assuming that $G$ is adjoint (resp. simply connected), see for example \cite{X4,X2}; this in turn gives rise to Springer correspondences for all $\Lg$ and $\Lg^*$ in any characteristic. 

If the characteristic $p$ of $\tk$ is very good for $G$ (namely $p$ is good, and $p$ does not divide $n+1$ if $G$ is of type $A_{n}$), there exists a $G$-equivariant isomorphism between $\cU_G$ and $\cN_\Lg$ by a theorem of Springer, moreover, there exists a $G$-invariant non-degenerate bilinear form on $\Lg$ which we can use to identify $\Lg$ and $\Lg^*$. Thus in such cases the sets $\cA_G$, $\cA_{\Lg}$ and $\cA_{\Lg^*}$ can be naturally identified. The Springer correspondence for $G$, i.e. the map $\gamma_G:\on{Irr}(W)\to\cA_G$ has been described explicitly in all characteristics \cite{L3,LS2,Sp,Sh,AL,S4}. When $G$ is of type $A_n$ and $p$ divides $n+1$, the maps $\gamma_G$, $\gamma_\Lg$, $\gamma_{\Lg^*}$ coincide with the corresponding maps for $GL(n+1)$.

Assume now that the characteristic  of $\tk$ is bad for $G$. In general we can no longer identify the sets $\cA_G$, $\cA_{\Lg}$ and $\cA_{\Lg^*}$. When $G$ is a classical group, the Springer correspondence for $\Lg$ and $\Lg^*$ in bad characteristic has been determined explicitly in  \cite{X3}. When $G$ is of type $F_4$ and $\on{char}\tk=2$, the Springer correspondence for $\Lg$ has been described by Spaltenstein \cite{S2}. When $G$ is of type $E_6,E_7$ or $E_8$ in bad characteristic, the Springer correspondence for $\Lg$ was partially described by the work of Spaltenstein \cite{S5} and Holt-Spaltenstein \cite{HS}. In this paper we describe the Springer correspondence  for $\Lg$ and $\Lg^*$ in the remaining cases when $G$ is of type $G_2$ or $F_4$. We will focus on the case of $\Lg^*$ with details given and treat the case of $\Lg$ briefly. The Springer correspondence for $\Lg$ and $\Lg^*$ turns out to be the same as in characteristic $0$ when $G$ is of type $G_2$ and $\on{char}\tk=2$, or when $G$ is of type $F_4$ and $\on{char}\tk=3$. 
The same holds for $\Lg^*$ when $G$ is of type $G_2$ and $\on{char}\tk=3$.

The paper is organized as follows. In \S\ref{sec-pre} we discuss some preliminaries. In particular, we determine the existence of  $G$-invariant non-degenerate bilinear forms on $\Lg$ (see Lemma \ref{bilinear}). In \S\ref{sec-nilp} we describe the component groups of centralizers for nilpotent elements in $\Lg^*$ and we compute the dimensions of the Springer fibers. In \S\ref{sec-sc} we describe the Springer correspondence  $\gamma_{\Lg^*}:\on{Irr}W\to\cA_{\Lg^*}$ for $\Lg^*$ explicitly. In particular, 
following \cite{Lu1}, we give an a priori description of the set $\gamma_{\Lg^*}^{-1}\{(\cO,\bar\bQ_l)\,|\,(\cO,\bar\bQ_l)\in\cA_{\Lg*}\}$, i.e. a set of Weyl group representations that parametrize nilpotent coadjoint orbits in $\Lg^*$. Finally, in \S\ref{sec-SCg} we describe the maps $\gamma_{\Lg}:\on{Irr}W\to\cA_\Lg$ for $\Lg$ explicitly.

{\bf Acknowledgement.} The author thanks the referee for carefully reading the manuscript and for many helpful comments.

\section{Preliminaries}\label{sec-pre}
In this section let $G$ be an exceptional group of type $G_2$ or $F_4$ defined over an algebraically closed field $\tk$.
\subsection{Structural constants and Bruhat decomposition} As in \cite{X3}, we make use of the generators and relations in $G$, $\Lg$, $\Lg^*$, and the Bruhat decomposition for $G$. We recall some  notations from \cite[\S7.3-\S7.5]{X3}, for details, see {\em loc.cit}. For $g\in G$ and $\xi\in\Lg^*$, we write $g.\xi$ for the coadjoint action, i.e., $g.\xi(x)=\xi(\on{Ad}(g^{-1})x)$ for all $x\in\Lg$, where $\on{Ad}$ denotes the adjoint action of $G$ on $\Lg$.

Fix a maximal torus $T$ of $G$ and $B\supset T$ a Borel subgroup of $G$. Let $R, R^+$ and $\Pi$ be the set of roots, positive roots, and simple roots, determined by $(G,T,B)$, respectively. For each $\alpha\in R$, let $x_\alpha:\mathbb{G}_a\to U_\alpha\subset G$ be an isomorphism  such that 
\beqn
\text{$sx_\alpha(t)s^{-1}=x_\alpha(\alpha(s)t)$ for all $s\in T$ and $t\in\mathbb{G}_a$,}
\eeqn  
\beqn
\text{ $n_\alpha(t):=x_\alpha(t)x_{-\alpha}(-t^{-1})x_\alpha(t)\in N_G(T)$ and $n_\alpha(1):=n_\alpha$ has image $w_\alpha$ in $W$,}
\eeqn
where $w_\alpha$ denotes the reflection with respect to $\alpha$. Define \beqn
 h_\alpha(t)=n_\alpha(t)n_\alpha(-1)\in T.
 \eeqn 
By Bruhat decomposition, each $g\in G$ can be written uniquely in the form 
\beq\label{bruhat}
\begin{gathered}
\text{$g=bn_wu_w$ for some  $b\in B$, ${u_w\in U_w:=\{\prod_{\alpha>0,w(\alpha)<0}x_\alpha(t_\alpha)\,|\,t_\alpha\in\mathbb{G}_a\}}$,}\\
\text{ and $n_w\in N_G(T)$ a representative of $w\in W$. }
\end{gathered}
\eeq

Let $\Lt$, $\Lb$, $\Ln$ be the Lie algebra of $T$, $B$, and the unipotent radical $U$ of $B$ respectively. Define $\Ln^*=\{\xi\in\Lg^*\,|\,\xi(\Lb)=0\}$. Then for any $\xi\in\N$, there exists $g\in G$ such that $g.\xi\in\Ln^*$. We make use of the Chevalley basis $\{h_\alpha,\,\alpha\in\Pi,\,e_\alpha,\,\alpha\in R\}$ of $\Lg$,  and the elements $\{e_\alpha'\in\Lg^*, \ \alpha\in R\}$ (in particular,  $\{e_\alpha,\,\alpha\in R^+\}$ is a basis of $\Ln$ and  $\{e_\alpha',\alpha\in R^+\}$ is a basis of $\Ln^*$), chosen in \cite[\S7.3]{X3},  where the adjoint action of $G$ is given as follows
\beq\label{commutator relations-l}
\begin{gathered}
\on{Ad}(x_\alpha(t))\,e_\beta=\sum_{i}\,t^i\,M_{\alpha,\beta,i}\,e_{i\alpha+\beta},\ \alpha,\beta\in R^+\\
\on{Ad}(h_\alpha(t))\,e_\beta=t^{A_{\alpha\beta}}e_\beta,\ \on{Ad}(n_\alpha)\,e_\beta=\eta_{\alpha,\beta}\, e_{w_\alpha(\beta)},\ \alpha\in\Pi,
\end{gathered}
\eeq
and the coadjoint action of $G$ is given as follows
\beq\label{commutator relations}
\begin{gathered}
x_\alpha(t).e_\beta'=\sum_{i}(-1)^i\,t^i\,M_{\alpha,-i\alpha-\beta,i}\,e_{i\alpha+\beta}',\ \alpha,\beta\in R^+\\
h_\alpha(t).e_\beta'=t^{A_{\alpha\beta}}e_\beta',\ n_\alpha.e_\beta'=\eta_{\alpha,\beta}\, e_{w_\alpha(\beta)}',\ \alpha\in\Pi.
\end{gathered}
\eeq
We will also make use of the following relations in $G$ (see for example \cite[\S 5.2, \S7.2]{C}), for $\alpha,\beta\in R$,
\beq\label{commutator relations-g}
\begin{gathered}
 x_\alpha(s)^{-1}x_\beta(t)^{-1}x_\alpha(s)x_\beta(t)=\prod_{i,j>0}x_{i\beta+j\alpha}(c_{ij\beta\alpha}\,(-t)^is^j),\  \alpha\neq\pm\beta\\
n_\alpha^2=h_\alpha(-1),\ n_\alpha\, x_{\beta}(t) n_\alpha^{-1}=x_{w_\alpha(\beta)}(\eta_{\alpha,\beta}\,t),
\ n_\alpha n_\beta n_\alpha^{-1}=h_{w_\alpha(\beta)}(\eta_{\alpha,\beta})\,n_{w_\alpha(\beta)}.
\end{gathered}
\eeq
Here 
$A_{\alpha,\beta}$, $M_{\alpha,\beta,i}$, $c_{ij\beta\alpha}$ and $\eta_{\alpha,\beta}(=\pm1)$ are structural constants (for their determination see \cite[\S7]{X3} and \cite[pages 77 and 94]{C}), in particular, they are integers.

\subsection{Invariant non-degenerate bilinear forms on $\Lg$}The following lemma was explained to the author by G. Lusztig. We include a proof here for completeness.
\begin{lemma}\label{bilinear}
Assume that $G$ is of type $G_2$ (resp. $F_4$). Then there exists a $G$-invariant non-degenerate bilinear form on $\Lg$ if and only if $\on{char}\tk\neq 3$ (resp. $\on{char}\tk\neq 2$).
\end{lemma}
Thus in the above cases, we have a $G$-equivariant isomorphism $\Lg\cong\Lg^*$. In turn the Springer correspondence maps $\gamma_\Lg$ and $\gamma_{\Lg^*}$ can be identified.

To prove Lemma \ref{bilinear}, we first note that the Lie algebra $\Lg$ is simple in the above cases \cite{H}.
 Assume that $G$ is of type $G_2$ and $\on{char}\tk\neq 3$. Let $\alpha$ be the short simple root and $\beta$ the long simple root. Define a symmetric bilinear form $(,)$ on $\Lg$ as follows
\beq\label{bilinear1}
\begin{gathered}
(h_\alpha,h_\alpha)=6,\ (h_\alpha,h_{\beta})=-3,\  (h_\beta,h_{\beta})=2\\
(e_\lambda,e_\mu)=0\text{ if }\lambda+\mu\neq0,\ (e_{\lambda},e_{-\lambda})=3\text{ (resp. 1) if $\lambda$ is a short (resp. long) root}.
\end{gathered}
\eeq 
Using \cite[Theorem 2.2]{Kac} one checks that \eqref{bilinear1} defines an invariant non-degenerate bilinear form on $\Lg$ over $\bC$, thus it is a scalar multiple of the Killing form, which is $G$-invariant. As all the structure constants are in $\bZ$, reducing mod $\on{char}\tk$, we see that the above bilinear form remains non-degenerate and $G$-invariant.

Assume that $G$ is of type $F_4$ and $\on{char}\tk\neq 2$.  Let $p,q$ be the long simple roots and $r,s$ the short simple roots such that $(q,r)\neq 0$. 
Define a symmetric bilinear form $(,)$ on $\Lg$ as follows
\beq\label{bilinear2}
\begin{gathered}
(h_p,h_p)=(h_q,h_q)=2,\ (h_p,h_q)=-1,\ (h_p,h_r)=(h_p,h_s)=(h_q,h_s)=0, \\  (h_q,h_r)=(h_r,h_s)=-2,\ (h_r,h_r)=(h_s,h_s)=4,\\ 
(e_\lambda,e_\mu)=0\text{ if }\lambda+\mu\neq0,\ (e_{\lambda},e_{-\lambda})=2\text{ (resp. 1) if $\lambda$ is a short (resp. long) root}.
\end{gathered}
\eeq 
The same argument as in the case of $G_2$ shows that the bilinear form $(,)$ of \eqref{bilinear2} is non-degenerate and $G$-invariant.

Now assume that $\on{char}\tk=3$ and $G$ is of type $G_2$, there are $6$ nilpotent orbits in $\Lg$ \cite{St} and $5$ nilpotent coadjoint orbits in $\Lg^*$ \cite{X1}.  Similarly, when $\on{char}\tk=2$ and $G$ is of type $F_4$, there are $22$ nilpotent orbits in $\Lg$ \cite{S2} and $18$ nilpotent coadjoint orbits in $\Lg^*$ \cite{X1}. Thus in these cases there can not exist a $G$-invariant non-degenerate bilinear form on $\Lg$.  This completes the proof of Lemma \ref{bilinear}.
\begin{remark}
Alternatively, as pointed out by the referee, it is known that $\Lg$ and $\Lg^*$ are not isomorphic as $G$-modules in these cases. Then Lemma \ref{bilinear} follows.
\end{remark}

\subsection{Weyl group representations and truncated induction}\label{sec W-rep}
Let $W$ be a Weyl group and let $\on{Irr}W$ denote the set of irreducible characters of $W$. For $\rho\in\on{Irr}W$,
 \beq\label{b function}
\begin{gathered}
\text{ let $b_\rho$ denote the minimal integer $d$ such that $\rho$ occurs in the $W$-module}
\\\text {$\mathfrak{S}_d(V)$ of all homogeneous polynomials of degree $d$ on the reflection space $V$.}
\end{gathered}
\eeq
For a subgroup $W'\subset W$, let $j_{W'}^W:\on{Irr}W'\to\on{Irr}W$ denote the truncated induction, i.e., for $\rho'\in\on{Irr}W'$, $j_{W'}^W\rho'$ is the unique $\rho\in\on{Irr}W$ such that $\rho$ occurs in $\on{Ind}_{W'}^W\rho'$, written as  $\langle\rho,\on{Ind}_{W'}^W\rho'\rangle\neq 0$, and $b_\rho=b_{\rho'}$.

 We label $\on{Irr}W$ using the notations of \cite{Al}. In particular, $\on{Irr}W(A_{n-1})=\on{Irr}S_n$ is labelled by partitions of $n$, where $(n)$ denotes the trivial character and $(1^n)$ denotes the sign character, and $\on{Irr}W(B_n)=\on{Irr}W(C_n)$ is labelled by a pair of partition $[\lambda:\mu]$. For $W$ of type $G_2$ or $F_4$, $\chi_{i,j}\in\on{Irr}W$ has degree $i$. For 
the reader's convenience, we list the values $b_\rho$ (see \eqref{b function}) for $\rho\in\on{Irr}W(G_2)$ and $\rho\in\on{Irr}W(F_4)$ in the following tables (see for example \cite{GP}).
\begin{table}[h]
\begin{center}
\bgroup
\def\arraystretch{1.2}
\begin{tabular}{>{$}l<{$} | >{$}l<{$} | >{$}l<{$} | >{$}l<{$} | >{$}l<{$} | >{$}l<{$} | >{$}l<{$} }
\hline
\rho & \chi_{1,1}& \chi_{2,1}&\chi_{2,2}& \chi_{1,3}&\chi_{1,4}&\chi_{1,2} \\\hline
b_\rho & 0&1&2&3&3&6 \\\hline
\end{tabular}
\egroup
\end{center}
\caption{$b$-invariants for $W(G_2)$}
\label{b inv G2}
\end{table}

\begin{table}[h]
\begin{center}
\bgroup
\def\arraystretch{1.2}
\begin{tabular}{ >{$}l<{$} | >{$}l<{$} | >{$}l<{$} | >{$}l<{$} | >{$}l<{$} | >{$}l<{$} | >{$}l<{$} | >{$}l<{$} | >{$}l<{$} | >{$}l<{$} | >{$}l<{$} | >{$}l<{$} | >{$}l<{$} | >{$}l<{$} }
\hline
\rho & \chi_{1,1}&\chi_{4,2}&\chi_{9,1}&\chi_{8,1}&\chi_{8,3}&\chi_{2,1}&\chi_{2,3}&\chi_{12,1}&\chi_{16,1}& \chi_{6,1}&\chi_{6,2}&\chi_{9,2}&\chi_{9,3} \\\hline
b_\rho & 0&1&2&3&3&4&4&4&5&6&6&6&6 \\\hline
\rho & \chi_{4,3}&\chi_{4,4}& \chi_{4,1}&\chi_{8,2}&\chi_{8,4}&\chi_{9,4}&\chi_{1,2}&\chi_{1,3}&\chi_{4,5}&\chi_{2,2}&\chi_{2,4}&\chi_{1,4} \\\hline
b_\rho & 7&7&8&9&9&10&12&12&13&16&16&24& \\\hline
\end{tabular}
\egroup
\end{center}
\caption{$b$-invariants for  $W(F_4)$}
\label{b inv F4}
\end{table}

\section{Nilpotent coadjoint orbits}\label{sec-nilp}
Let $G$ continue to be an exceptional group of type $G_2$ or $F_4$ over an algebraically closed field $\tk$. In this section, we describe the component groups $A_G(\xi):=Z_G(\xi)/Z_G(\xi)^0$ for $\xi\in\N$.   Moreover, we compute the dimensions of Springer fibers. More precisely, let $B$ be a fixed Borel subgroup of $G$ with unipotent radical $U$. We write $\Lb=\on{Lie}B$ and $\Ln=\on{Lie}U$. Define
$$\Lb^*=\{\zeta\in\Lg^*\,|\,\zeta(\Ln)=0\}\text{ and }\cB^G_\xi=\{gB\in G/B\,|\,g^{-1}.\xi\in\Lb^*\},\ \xi\in\Lg^*.$$
We show in \S\ref{ssec-sf} that 
\begin{proposition}\label{dimension}
For $\xi\in\N$,
\beq\label{Sfiber}
\dim\cB^G_\xi=\frac{\dim Z_G(\xi)-\on{rank}G}{2}.
\eeq
\end{proposition}

The component groups $A_G(x)$ for $x\in\cN_\Lg$ have been determined in all characteristics, as well as the dimensions of Springer fibers, see \S\ref{ssec-nilp orbit} for exposition. In view of Lemma \ref{bilinear}, it suffices to consider the cases when $G$ is an exceptional group of type $G_2$ (resp. $F_4$) over an algebraically closed field $\tk$ of characteristic $3$ (resp. $2$), which we will assume in the remainder of this section. 

The nilpotent coadjoint orbits in $\Lg^*$ have been classified in \cite[\S7]{X1}.  We use the same notation and representatives for the orbits as in {\em loc.cit.} In particular, the orbits are labeled as in the Bala-Carter classification (while the two  new orbits in bad characteristic are labeled as $(B_3)_2$ and $(A_2)_2$). The tilde in $\widetilde{A}_1$ etc.
indicates short roots.

The results are listed in Table \ref{tab:G2 coadjoint} when $G$ is of type $G_2$ and $\on{char}\tk=3$, and in Table \ref{tab:F4 coadjoint} when $G$ is of type $F_4$ and $\on{char}\tk=2$.
\begin{table}[h]
\centering
\bgroup
\def\arraystretch{1.2}
\begin{tabular}{c|c|c|c|c }
\hline
Orbit&representative $\xi$&$\dim Z_G(\xi)$&$A_G(\xi)$&$\dim\cB_\xi$\\
\hline
$G_2$&$e_\alpha'+e_\beta'$&$2$&1&0\\
$G_2(a_1)$&$e_\beta'+e_{2\alpha+\beta}'$&$4$&$S_3$&1\\
$\widetilde{A_1}$&$e_\alpha'$&$6$&1&2\\
${A_1}$&$e_\beta'$&$8$&1&3\\
$\emptyset$ &0&$14$&1&6\\
\hline
\end{tabular}
\egroup
\vspace{.1in}
\caption{Type $G_2$, nilpotent coadjoint orbits, $\on{char}\tk=3$}
\label{tab:G2 coadjoint}
\end{table}

\begin{table}[h]
\centering
\bgroup
\def\arraystretch{1.2}
\begin{tabular}{c|c|c|c|c}
\hline
Orbit&representative $\xi$&$\dim Z_G(\xi)$&$A_G(\xi)$&$\dim \cB_\xi$\\
\hline
$F_4$&$e_p'+e_q'+e_r'+e_s'$&$4$&1&0\\
$F_4(a_1)$&$e_p'+e_{qr}'+e_{q2r}'+e_s'$&$6$&$S_2$&1\\
$F_4(a_2)$&$e_p'+e_{qr}'+e_{rs}'+e_{q2r2s}'$&$8$&$1$&2\\
$B_3$&$e_p'+e_{qrs}'+e_{q2r}'+e_{pq2rs}'$&${10}$&1&3\\
$C_3$&$e_s'+e_{q2r}'+e_{pqr}'$&$10$&1&3\\
$F_4(a_3)$&$e_{pqr}'+e_{qrs}'+e_{pq2r}'+e_{q2r2s}'$&${12}$&$S_4$&4\\
$(B_3)_2$&$e_p'+e_{qr}'+e_{q2r2s}'$&$12$&1&4\\
$C_3(a_1)$&$e_{pqr}'+e_{q2rs}'+e_{q2r2s}'$&$14$&$S_2$&5\\
$B_2$&$e_{pqr}'+e_{q2r2s}'$&$16$&$S_2$&6\\
$\widetilde{A_2}+A_1$&$e_{pqr}'+e_{q2rs}'+e_{p2q2r2s}'$&$16$&1&6\\
$A_2+\widetilde{A_1}$&$e_{p2q2r}'+e_{q2r2s}'+e_{pq2rs}'$&$18$&1&7\\
$\widetilde{A_2}$&$e_{pqrs}'+e_{q2rs}'$&$22$&1&9\\
$A_2$&$e_{p2q2r}'+e_{pq2r2s}'+e_{p2q3r2s}'$&$22$&1&9\\
$A_1+\widetilde{A_1}$&$e_{p2q2r2s}'+e_{p2q3rs}'$&$24$&1&10\\
$(A_2)_2$&$e_{p2q2r}'+e_{pq2r2s}'$&$28$&1&12\\
$\widetilde{A_1}$&$e_{p2q3r2s}'$&$30$&$S_2$&13\\
$A_1$&$e_{2p3q4r2s}'$&$36$&1&16\\
$\emptyset$&$0$&$52$&1&24\\
\hline
\end{tabular}
\egroup
\vspace{.1in}
\caption{Type $F_4$, nilpotent coadjoint orbits, $\on{char}\tk=2$}
\label{tab:F4 coadjoint}
\end{table}

\subsection{Component groups of centralizers} To describe the structure of the component groups $A_G(\xi)$, $\xi\in\N$,  we make use of the following lemma.
\begin{lemma}[{\cite[3.4, p176]{B}}]Let $\sigma:G\to G$ be the Frobenius endomorphism induced by $\tk\to\tk:a\mapsto a^q$. Let $\cO\subset\Lg^*$ be a $G$-orbit. Then the set of $G^\sigma$-orbits in $\cO^\sigma$ is in bijection with $H^1(\sigma,A_G(\xi))$, where $\xi\in\cO$. 
\end{lemma}
Note that for $q$ large enough, $\sigma$ acts trivially on $A_G(\xi)$. Thus $|H^1(\sigma,A_G(\xi))|$ equals the number of conjugacy classes in $A_G(\xi)$. It then follows from \cite[\S7, Table 1 and Table 2]{X1} that we only need to consider the  representative $\xi=e_\beta'+e_{2\alpha+\beta}'$ for the orbit $G_2(a_1)$ and the representative $\xi_5=e_{pqr}'+e_{qrs}'+e_{pq2r}'+e_{q2r2s}'$ for the orbit $F_4(a_3)$. Here we have used the fact that $S_2\cong\bZ/2\bZ$ is the only finite group with exactly 2 conjugacy classes.

Assume that $G$ is of type $G_2$ and $\on{char}\tk=3$. We determine $A_G(\xi)$ for $\xi=e_\beta'+e_{2\alpha+\beta}'$. Using \eqref{bruhat} and \eqref{commutator relations}, one can check that $Z_G(\xi)\subset B$ and 
\beq\label{component group 1}
A_G(\xi)\cong\langle\gamma_1,\gamma_2\rangle\cong S_3,
\eeq where $\gamma_1=h_\beta(-1)$ and $\gamma_2=h_\beta(-1)x_\alpha(\eta)$ ($\eta^2=-1$). In particular, we have used the following formula
\beqn
u(\mathbf{t}).\xi=\xi+(t_1^3+t_1)e'_{3\alpha+\beta}+(t_5+t_3)e'_{3\alpha+2\beta}.
\eeqn
where
$
u(\mathbf{t})=x_\alpha(t_1)x_\beta(t_2)x_{\alpha+\beta}(t_3)x_{2\alpha+\beta}(t_4)x_{3\alpha+\beta}(t_5)x_{3\alpha+2\beta}(t_6),
$ $t_i\in\mathbb{G}_a$.

Assume that $G$ is of type $F_4$ and $\on{char}\tk=2$. We determine $A_G(\xi_5)$ for $\xi_5=e_{pqr}'+e_{qrs}'+e_{pq2r}'+e_{q2r2s}'$. Let $P$ be the standard parabolic group $P_J$, where $J=\{p,r,s\}$. Then using \cite[Theorem 7.3]{CP} and the description of nilpotent pieces in $\Lg^*$ in \cite[\S7.6]{X1} we see that $Z_G(\xi_5)\subset P$. Using \eqref{bruhat} and \eqref{commutator relations}, one checks that
\beq\label{component group s4}
Z_G^0(\xi_5)=Z_{U_P}(\xi_5),\ \ A_G(\xi_5)\cong\langle\gamma_1,\gamma_2,\gamma_3\rangle\cong S_4,
\eeq
where
\beqn
\gamma_1=x_p(1)\,x_s(1),\ \ \gamma_2=n_p\,n_s,\ \ \gamma_3=x_p(1)\,x_s(1)\,x_r(1)\,n_r\,x_r(1)
\eeqn
with
\beqn
\gamma_1^2=\gamma_2^2=\gamma_3^2=1,\ (\gamma_1\gamma_2)^3=(\gamma_2\gamma_3)^3=1,\ (\gamma_1\gamma_3)^2=1.
\eeqn
In particular we have used the following formula (here the assumption $\on{char}\tk=2$ is used)
\beqn
\begin{gathered}
u(\mathbf{t}).\xi_5=e_{pqr}'+e_{qrs}'+(1+t_2)e_{pq2r}'+t_2\,e_{q2rs}'+(t_1+t_3)\,e_{pqrs}'+(1+t_4)e_{q2r2s}'\\
+(t_1t_2+t_2t_3+t_4)e_{pq2rs}'+(t_1+t_3^2+t_4(t_1+t_3))e_{pq2r2s}'
\end{gathered}
\eeqn
where $u(\mathbf{t})=x_p(t_1)\,x_r(t_2)\,x_s(t_3)\,x_{rs}(t_4)$, $t_i\in\mathbb{G}_a$.
\begin{remark}
For the element $\xi_5$ above, we have another interesting way to deduce that $A_G(\xi_5)\cong S_4$ as follows. Taking $q$ big enough, we can choose a representative $\xi$ of the orbit of $\xi_5$ such that the Frobenius morphism $\sigma$ fixes $\xi$ and acts trivially on $A_G(\xi)$. This implies that $|Z_G(\xi)^\sigma|=|A_G(\xi)||(Z_G(\xi)^0)^\sigma|$. In view of \cite[Table 2]{X1}, we conclude that $A_G(\xi)$ is a finite group of order $24$ with exactly $5$ conjugacy classes. It turns out that the finite groups with exactly $5$ conjugacy classes have been classified by Miller \cite{M} in 1919; the list is the following: the alternating group $A_5$, $S_4$, the non-cyclic group of order $21$, the metacylic group of order $20$, the dihedral group of order $14$, the dihedral group of order $8$, the quaternion group, and $\bZ/5\bZ$. We can then easily  identify our group from the list.
\end{remark}

\subsection{Springer fibers} \label{ssec-sf}
 Let $L$ denote a Levi subgroup of a proper parabolic subgroup $P$ of $G$. We write $\Ll=\on{Lie}L$, $\Ln_P=\on{Lie}U_P$, where $U_P$ is the unipotent radical of $P$. Let 
 \beqn
\text{ $\Ll^*=\{\zeta\in\Lg^*\,|\,\zeta(\Ln_P\oplus\Ln_P^-)=0\}$ and $\Ln_P^*=\{\zeta\in\Lg^*\,|\,\zeta(\Ll\oplus\Ln_P)=0\}$,}
 \eeqn
 where $\Lg=\Ll\oplus\Ln_P\oplus\Ln_P^-$. 
 
 Let $\xi\in\N$. We prove Proposition \ref{dimension} using the following statements. 
\beq\label{red-a}
\begin{gathered}
\text{ Assume that $\xi$ lies in $\Ll^*$ and $\dim Z_L(\xi)=2\dim\cB_\xi^L+\on{rank} L$.}\\\text{ Then $\dim Z_G(\xi)=2\dim\cB_\xi^G+\on{rank}G$.} 
\end{gathered}
\eeq
\beq\label{red-b}
\begin{gathered}
\text{ If $\xi$ lies in the orbit obtained by inducing the orbit of $\xi'\in\Ll^*$ \footnotemark, and }\\\text{ $\dim Z_L(\xi')=2\dim\cB_{\xi'}^L+\on{rank}L$, then $\dim Z_G(\xi)=2\dim\cB_\xi^G+\on{rank}G$.} 
 \end{gathered}
 \footnotetext{see \cite[\S4]{X1}, i.e., $\cO_{\xi}\cap(\cO_{\xi'}^L+\Ln_P^*)$ is dense in $\cO_{\xi'}^L+\Ln_P^*$, where $\cO_{\xi'}^L$ denotes the $L$-orbit of $\xi'$; we write $\cO_\xi=\on{Ind}_{\Ll^*}^{\Lg^*}\cO_{\xi'}^L$.}
  \eeq
The statement~\eqref{red-a}  is proved using the same argument as in \cite[II 3.14]{Spa1}. In particular, one uses that the equality in \eqref{Sfiber} holds if there exists $w\in W$ such that $\cO_\xi\cap\Ln^*\cap\Ln_w^*$ is dense in $\Ln^*\cap\Ln_w^*$, where $\Ln_w^*=\{\zeta\in\Lg^*\,|\,\zeta(w\Lb w^{-1})=0\}$. The statement~\eqref{red-b} is proved using the same argument as in \cite[Theorem 1.3]{LS}. In particular,  the proof shows if $\cO_\xi=\on{Ind}_{\Ll^*}^{\Lg^*}\cO_{\xi'}^L$, then
\beq\label{dim-induced}
\dim\cO_\xi=\dim\cO_{\xi'}^L+2\dim U_P.
\eeq

Note that  by \cite[Propositon 3.1]{X1}, the equality \eqref{Sfiber} holds  for any $\xi'\in\cN_{\Ll^*}$ (replacing $G$ by $L$).
Thus it suffices to show that 
\beq\label{orbits}
\begin{gathered}
\text{
every $\xi\in\N$ either lies in $\Ll^*$  or lies in an orbit induced from a }\\\text{nilpotent orbit in $\Ll^*$ for a proper Levi subgroup $L$. }
\end{gathered}
\eeq
For $G_2$, one readily checks~\eqref{orbits}. For example, the orbit $G_2$ is induced from the zero orbit in a maximal torus, the orbit $G_2(a_1)$ is induced from the zero orbit in a Levi subgroup of type $\widetilde{A}_1$. For $F_4$, we make use of the following table for induced orbits. The first row contains orbits $\widetilde\cO$ in $\N$, the second row indicates the type of the Levi subgroup $L$, and the third row contains orbits $\cO$ in $\cN_{\Ll^*}$ such that $\widetilde{\cO}=\on{Ind}_{\Ll^*}^{\Lg^*}\cO$, where $0$ denotes the zero orbit and  we use the notation in \cite{X2} for orbits in a Levi of type $B_3$ or $C_3$. The last row of the table is included for use in \S\ref{ssec-F4} and will be explained there. We remark that we did not attempt to identify all induced orbits here, rather only the cases that are needed for our purpose.
\begin{table}[h]
\begin{center}
\bgroup
\def\arraystretch{1.3}
\begin{tabular}{ >{$}l<{$} | >{$}l<{$} | >{$}l<{$} | >{$}l<{$} | >{$}l<{$} | >{$}l<{$} | >{$}l<{$} | >{$}l<{$} | >{$}l<{$}  }
\widetilde\cO&\widetilde{A}_2&A_2&B_2&F_4(a_3)&F_4(a_3)&B_3&F_4(a_1)&F_4(a_2)\\\hline
\text{type of }L&B_3&C_3&C_3&\widetilde{A}_2A_1&B_3&B_3&B_3&A_1\widetilde{A}_1\\\hline
\cO&0&0&1_1^6&0&(1;1^4_1)&(0;3^2_3)&(2;1^2_1)&0\\
\hline
\rho_{\xi',1}^L&[-:1^3]&[-:1^3]&[1^3:-]&[1^3]\times[1^2]&[1:1^2]&[-:3]&[2:1]&[2]\times[2]
\end{tabular}
\egroup
\end{center}
\vspace{.1in}
\caption{Induced nilpotent coadjoint orbits, type $F_4$, $\on{char}\tk=2$}
\label{induced}
\end{table}
As an example, we explain the case of the orbit $F_4(a_3)$ being induced from the orbit $(1;1^4_1)$ in $\Ll^*$ for $L$ of type $B_3$; the other cases can be checked similarly. Let $\xi=e_{pqr}'+e_{qrs}'+e_{pq2r}'+e_{q2r2s}$, which is in the orbit $F_4(a_3)$, and let $\xi'=e_{pqr}'+e_{pq2r}'$. Let $P$ be the standard parabolic subgroup $P_J$, $J=\{p,q,r\}$ and $L$ the standard Levi subgroup of $P$. Then $L$ is of type $B_3$. Moreover, $\xi'\in\Ll^*$ and $\xi\in\xi'+\Ln_P^*$. Using \cite{X2}, one checks that $\xi'$ lies in the orbit $(1;1^4_1)$ of $\Ll^*$, moreover,  $\dim\cO_{\xi'}^L=10$. Note that $\dim\cO_\xi=40=\dim\cO_{\xi'}+2\dim U_P=\dim\on{Ind}_{\Ll^*}^{\Lg^*}\cO_{\xi'}^L$ (in the last equality we use \eqref{dim-induced}). Since $\on{Ind}_{\Ll^*}^{\Lg^*}\cO_{\xi'}^L$ is the unique dense orbit in $G.(\cO_{\xi'}+\Ln_P^*)$, we obtain that 
\beq\label{ex-ind}
\cO_\xi=\on{Ind}_{\Ll^*}^{\Lg^*}\cO_{\xi'}^L,\ \xi=e_{pqr}'+e_{qrs}'+e_{pq2r}'+e_{q2r2s},\ \xi'=e_{pqr}'+e_{pq2r}'.
\eeq 
Now one checks readily that~\eqref{orbits} holds for $F_4$. This completes the proof of Proposition \ref{dimension}.

\section{The Springer correspondence for $\Lg^*$}\label{sec-sc}
Let $G$ be a connected simply connected algebraic group defined over an algebraically closed field $\tk$ of prime characteristic. Fix a prime $l\neq\on{char}\tk$. Let 
\beq\label{Smap}
\gamma_{\Lg^*}:\on{Irr}W\to\cA_{\Lg^*}
\eeq
be the injective Springer correspondence map constructed for $\Lg^*$ as in \cite[\S5]{X2}, where $\on{Irr}W$ and $\cA_{\Lg^*}$ are the sets defined in the introduction. Namely, $\on{Irr}W$ denotes the set of irreducible characters of the Weyl group $W$ of $G$ and $\cA_{\Lg^*}$ is the set of all pairs $(\cO,\cE)$ where $\cO\subset\N$ is a $G$-orbit and $\cE$ is an irreducible $G$-equivariant $\bar\bQ_l$-local system on $\cO$ (up to isomorphism).  The constructions and proofs in \cite[\S5]{X2} apply for any simply connected $G$. In what follows we briefly recall the construction. Consider the following proper maps
\beqn
\varphi_0:G\times^B\Ln^*\to\N,\ \ \varphi:G\times^B\Lb^*\to\Lg^*.
\eeqn
The map $\varphi_0$ is semismall and the map $\varphi$ is  small. Moreover,
\beq\label{decomposition}
\varphi_{0*}\bar\bQ_l[-]\cong\bigoplus_{(\cO,\cE)\in\cA_{\Lg^*}}\on{IC}(\bar\cO,\cE)\otimes V_{\cO,\cE}\cong\varphi_*\bar\bQ_l|_{\N}[-]\cong\bigoplus_{\rho\in\on{Irr}W}\rho\otimes(\varphi_*\bar\bQ_l)_\rho|_{\N}[-]
\eeq
where $[-]$ denotes shift by $\dim\N$, $\on{IC}(\bar\cO,\cE)$ is the perverse IC-extension of the local system $\cE$ on $\cO$, and $(\varphi_*\bar\bQ_l)_\rho=\on{Hom}_{\bar\bQ_l[W]}(\rho,\varphi_*\bar\bQ_l)$. The map $\gamma_{\Lg^*}$ maps $\rho\in\on{Irr}W$ to the unique pair $(\cO,\cE)\in\cA_{\Lg^*}$ such that $\on{IC}(\bar\cO,\cE)\cong(\varphi_*\bar\bQ_l)_\rho|_{\N}[-]$. Note that this implies that $\rho\cong V_{\cO,\cE}$ under $W$-action.

In this section we describe the map $\gamma_{\Lg^*}$ of \eqref{Smap} explicitly assuming that $G$ is of type $G_2$ (resp. $F_4$) and $\on{char}\tk=3$ (resp. $\on{char}\tk=2$).  The results are given in Table \ref{tab:G2SCg*} (resp. Table \ref{tab:F4SCg*}). For a pair $(\cO,\cE)\in\cA_{\Lg^*}$, we write $\rho_{\xi,\phi}^G$ (or simply $\rho_{\xi,\phi}$) for the inverse image $\gamma_{\Lg^*}^{-1}(\cO,\cE)$, where $\xi\in\cO$ and $\phi\in \on{Irr}A_G(\xi)$ corresponds to the local system $\cE$.  In particular, as mentioned in the introduction, the Springer correspondence for $\Lg^*$ turns out to be the same as in characteristic $0$ when $G$ is of type $G_2$, or when $G$ is of type $F_4$ and $\on{char}\tk=3$. 

\begin{table}[h]
\centering
\bgroup
\def\arraystretch{1.3}
\begin{tabular}{c c c c c c c}
\hline
Orbit of $\xi$&$\phi$&$\rho_{\xi,\phi}$& \vline\ \vline& Orbit of $\xi$ &$\phi$ &$\rho_{\xi,\phi}$\\
\hline
$G_2$&$(1)$&$\chi_{1,1}$& \vline\ \vline&$\widetilde{A}_1$&$(1)$&$\chi_{2,2}$\\
$G_2(a_1)$&$(3)$&$\chi_{2,1}$& \vline\ \vline&${A_1}$&$(1)$&$\chi_{1,4}$\\
$G_2(a_1)$&$(2,1)$&$\chi_{1,3}$&\vline\ \vline&
$\emptyset$ &(1)&$\chi_{1,2}$ \\
$G_2(a_1)$&$(1^3)$&$-$&\vline\ \vline&\\
\hline
\end{tabular}
\egroup
\vspace{.1in}
\caption{Springer correspondence for $\Lg^*$, type $G_2$ and $\on{char}\tk=3$}
\label{tab:G2SCg*}
\end{table}
\begin{table}[h]
\centering
\bgroup
\def\arraystretch{1.3}
\begin{tabular}{c c c c c c c c ccccccc}
\hline
Orbit of $\xi$ &$\phi$ &$\rho_{\xi,\phi}$& \vline\ \vline& Orbit of $\xi$ &$\phi$ &$\rho_{\xi,\phi}$& \vline\ \vline& Orbit of $\xi$ &$\phi$ &$\rho_{\xi,\phi}$\\
\hline
$F_4$&$(1)$&$\chi_{1,1}$&    \vline\ \vline&$F_4(a_3)$&$(2,1,1)$& ${\chi_{1,3}}$&    \vline\ \vline&$\widetilde{A_2}$&$(1)$&$\chi_{8,2}$\\
$F_4(a_1)$&$(2)$&$\chi_{4,2}$& \vline\ \vline&$F_4(a_3)$&$(1^4)$& $-$&    \vline\ \vline&$A_2$&$(1)$&$\chi_{8,4}$\\
$F_4(a_1)$&$(1^2)$ &$\chi_{2,3}$& \vline\ \vline&$(B_3)_2$&$(1)$&$\chi_{2,1}$&    \vline\ \vline&$A_1+\widetilde{A_1}$&$(1)$&$\chi_{9,4}$\\
$F_4(a_2)$&$(1)$&$\chi_{9,1}$& \vline\ \vline&$C_3(a_1)$&$(2)$&$\chi_{16,1}$&    \vline\ \vline&$(A_2)_2$&$(1)$&$\chi_{1,2}$\\
$B_3$&$(1)$&$\chi_{8,1}$& \vline\ \vline&$C_3(a_1)$&$(1^2)$ &${\chi_{4,4}}$&    \vline\ \vline&$\widetilde{A_1}$&$(2)$&$\chi_{4,5}$\\
$C_3$&$(1)$&$\chi_{8,3}$& \vline\ \vline&$B_2$&$(2)$&$\chi_{9,2}$&    \vline\ \vline&$\widetilde{A_1}$&$(1^2)$ &$\chi_{2,2}$\\
$F_4(a_3)$&$(4)$&$\chi_{12,1}$& \vline\ \vline&$B_2$&$(1^2)$ &$\chi_{4,1}$&    \vline\ \vline&$A_1$&$(1)$&$\chi_{2,4}$\\
$F_4(a_3)$&$(3,1)$& ${\chi_{9,3}}$& \vline\ \vline&$\widetilde{A_2}+A_1$&$(1)$&$\chi_{6,1}$&    \vline\ \vline&$\emptyset$&$(1)$&$\chi_{1,4}$\\
$F_4(a_3)$&$(2,2)$& ${\chi_{6,2}}$& \vline\ \vline&$A_2+\widetilde{A_1}$&$(1)$&$\chi_{4,3}$&    \vline\ \vline&\\
\hline
\vspace{.1in}
\end{tabular}
\egroup
\caption{Springer correspondence for $\Lg^*$, type $F_4$ and $\on{char}\tk=2$}
\label{tab:F4SCg*}
\end{table}

We remark that in view of the results in \cite{X3}, Tables \ref{tab:G2SCg*} and  \ref{tab:F4SCg*}, the pair $(G_2(a_1),\cE_{(1^3)})$ for $G_2$ and the pair $(F_4(a_3),\cE_{(1^4)})$ for $F_4$ are the cuspidal pairs in the sense of \cite{L3}, where $\cE_\phi$ denotes the local system corresponding to $\phi\in\on{Irr}A_G(\xi)$.

\subsection{The methods}

We describe the map $\gamma_{\Lg^*}$ of \eqref{Smap} following the methods in \cite{AL,S4,S3}. In particular, we use extensively the tables of induce/restrict matrix for Weyl group representations given by Alvis in \cite{Al}. More specifically,  we make use of  the following statements, for which the proofs follow the same arguments as in \cite{AL,S4,S3} and the references there. Thus we only briefly comment on the proofs. Throughout this subsection $\xi\in\N$, $L$ is a Levi subgroup  of a proper parabolic subgroup $P$ of $G$, $\Ll=\on{Lie}L$, $U_P$ is the unipotent radical of $P$, and $W_L$ is the Weyl group of $L$,  regarded naturally as a subgroup of $W$. 

To begin with, let $b_\rho$ be defined as in \eqref{b function} for $\rho\in\on{Irr}W$. Using the same argument as in, for example, \cite[\S13.12]{J}, i.e.  taking stalks at 0 in \eqref{decomposition}, using Proposition \ref{dimension} and the fact that $[\rho:H^{2i}(\cB,\bar\bQ_l)]=[\rho:\mathfrak{S}_i(V)]$ (this follows from the facts that $H^{\bullet}(\cB,\bar\bQ_l)$ is isomorphic to the coinvariant algebra $C$ of $\mathfrak{S}(V)$ and that $C\otimes\mathfrak{S}(V)^W\cong \mathfrak{S}(V)$), we see that 
 \beq\label{b function and Springer fiber}
\text{ $b_{\rho_{\xi,1}^{G}}=\dim\cB_\xi^G$.}
\eeq
\indent Next let $\xi'\in\cN_{\Ll^*}$. Consider the permutation representation $\varepsilon_{\xi,\xi'}$ of $A_G(\xi)\times A_L(\xi')$ afforded by the finite set $S_{\xi,\xi'}$ of irreducible components  of dimension $d_{\xi,\xi'}$ of $Y_{\xi,\xi'}:=\{g\in G\,|\,g^{-1}.\xi\in\xi'+\Ln_P^*\}$ \footnote{Here $d_{\xi,\xi'}=(\dim Z_G(\xi)+\dim Z_L(\xi'))/2+\dim U_P$.}. The same proof as in \cite{L3} (see also \cite{X3}) shows that for $\phi\in \on{Irr}A_G(\xi)$ and $\phi'\in\on{Irr} A_L(\xi')$, we have
\beq\label{restriction formula}
\langle \phi\otimes\phi',\varepsilon_{\xi,\xi'}\rangle=\langle{\on{Res}_{W_L}^W\rho_{\xi,\phi}^G,\rho_{\xi',\phi'}^L}\rangle_{W_L}:=n_{\xi,\xi',\phi,\phi'}.
\eeq
In particular, we make use of the following special cases of \eqref{restriction formula}.
\beq\label{restriction}
\text{ If $n_{\xi,\xi',1,1}=0$, then $S_{\xi,\xi'}=\emptyset$ and all $n_{\xi,\xi',\phi,\phi'}=0$.}
\eeq

If the orbit of $\xi$ is induced from that of $\xi'\in\Ll^*$, then we can assume that $\xi\in\xi'+\Ln_P^*$. Using the same argument as in \cite[1.3, 1.5]{LS} and \cite{S3}, we see that in this case 
\beq\label{perm rep}
\text{$S_{\xi,\xi'}$ is isomorphic to $A_G(\xi)/N$ as sets with $A_G(\xi)\times A_L(\xi')$-actions,}
\eeq
where 
\beq\label{groups N H}
\text{$N=Z_{Z_L^0(\xi')U_P}(\xi)/Z_G^0(\xi)$ is a normal subgroup of $H=Z_P(\xi)/Z_G^0(\xi)$,  $H/N\cong A_L(\xi')$,}
\eeq
and the group $A_G(\xi)\times H/N\cong A_G(\xi)\times A_L(\xi')$ acts on $A_G(\xi)/N$ by $(a,hN).(xN)=axh^{-1}N$.
Thus it follows from \eqref{restriction formula} and \eqref{perm rep} that for $\phi\in \on{Irr}A_G(\xi)$ and $\phi'\in\on{Irr} A_L(\xi')$, 
\beq\label{restriction of component gp}
\text{$n_{\xi,\xi',\phi,\phi'}=\langle\tilde{\phi'},\on{Res}^{A_G(\xi)}_H\phi\rangle$, where $\tilde{\phi'}$ is the lift of $\phi'$ to $H$.}
\eeq
Moreover, the same argument as in \cite[1.3]{LS} shows that $\dim \cB_\xi^G=\dim\cB_{\xi'}^L$. It then follows from \eqref{b function and Springer fiber}, the definition of the truncated induction operator $j$ (see \S\ref{sec W-rep}) and the above discussion that
\beq\label{induction}
\text{ If the orbit of $\xi$ is induced from that of $\xi'\in\cN_{\Ll^*}$,  then  $\rho_{\xi,1}^G=j_{W_L}^W(\rho_{\xi',1}^L)$.}
 \eeq

Finally, the same argument as in \cite{AL} shows that
\beq\label{Levi}
\begin{gathered}
\text{If $\xi\in\Ll^*$, then $\langle\rho_{\xi,1}^G,\on{Ind}_{W_L}^W\hat{\rho}_{\xi,1}^L\rangle\neq 0$,}\text{ where $\hat{\rho}_{\xi,1}^L=\sum(-1)^iH^i(\cB^L_\xi,\bar\bQ_l)$}.
\end{gathered}
\eeq

In addition, we remark that the decomposition \eqref{decomposition} implies that we have a decomposition 
\beqn
H^{2d_\xi}(\cB_\xi^G,\bar\bQ_l)\cong\bigoplus_{\phi\in \on{Irr}A_G(\xi)}\phi\otimes \rho^G_{\xi,\phi}
\eeqn
as $A_G(\xi)\times W$ representations, where $d_\xi=\dim\cB_\xi^G$ and $A_G(\xi)$ acts via the permutation representation $\epsilon_\xi$ afforded by the set of irreducible components of $\cB_\xi^G$. In particular, we see that
\beq\label{perm-rep}
\rho_{\xi,\phi}^G\neq 0\text{ if and only if }\phi\text{ occurs in the above permutation representation $\epsilon_\xi$.}
\eeq

\subsection{Type $G_2$ in characteristic 3}We verify in this subsection the description of the map $\gamma_{\Lg^*}$ given in Table \ref{tab:G2SCg*}.  First, using \eqref{b function and Springer fiber}, Table \ref{b inv G2} and Table \ref{tab:G2 coadjoint}, one determines easily the images of $\chi_{1,1}$, $\chi_{1,2}$, $\chi_{2,1}$ and $\chi_{2,2}$ under $\gamma_{\Lg^*}$.

For $\xi_4=e_\beta'$ in the orbit $A_1$, note that $\xi_4$ is a regular nilpotent element in $\Ll^*$, where $L$ is a Levi subgroup of type $A_1$.   Using \eqref{Levi}, we see that $\langle\rho_{\xi_4,1},\on{Ind}_{W(A_1)}^W[2]\rangle\neq 0$. Thus it follows from \cite[Table 64]{Al} that $\rho_{\xi_4,1}=\chi_{1,4}$. 

It remains to determine $\gamma_{\Lg^*}(\chi_{1,3})$. To this end, it is enough to determine the pair $(\xi_2,\phi)$, $\xi_2=e_\beta'+e_{2\alpha+\beta}'$ in the orbit $G_2(a_1)$ and $\phi\neq 1$, that appears in the image of $\gamma_{\Lg^*}$. In view of \eqref{perm-rep} we study the permutation representation $\epsilon_{\xi_2}$ of $A_G(\xi_2)\cong S_3$ afforded by the set of irreducible components of $\cB_{\xi_2}$. Using the Bruhat decomposition \eqref{bruhat} and \eqref{commutator relations}, one can check that $g^{-1}.\xi_2\in\Lb^*$ if and only if $g^{-1}\in B$, $g^{-1}\in Bs_\alpha x_\alpha(t_1)$, or $g^{-1}\in Bs_\beta s_\alpha x_\alpha(\varpi)x_{3\alpha+\beta}(t_5)$, where $\varpi^3+\varpi=0$. Thus $\cB_{\xi_2}$ has $4$ irreducible components. Using the description of $A_G(\xi_2)$ in \eqref{component group 1}, one readily checks that the permutation representation $\epsilon_{\xi_2}$ decomposes as $(3)\oplus(3)\oplus (2,1)$. Thus $\gamma_{\Lg^*}(\chi_{1,3})=(\xi_2,(2,1))$. This completes the verification of Table \ref{tab:G2SCg*}.

\subsection{Type $F_4$ in characteristic $2$} \label{ssec-F4}

We verify in this subsection the description of the map $\gamma_{\Lg^*}$ given in Table \ref{tab:F4SCg*}.  First, using \eqref{b function and Springer fiber}, Table \ref{b inv F4} and Table \ref{tab:F4 coadjoint}, one determines easily the image of $\chi_{1,1}$, $\chi_{4,2}$, $\chi_{9,1}$, $\chi_{16,1}$,  $\chi_{9,4}$,  $\chi_{4,5}$, and  $\chi_{1,4}$, under $\gamma_{\Lg^*}$.

Using \eqref{b function and Springer fiber}, \eqref{induction}, \cite[Tables 60-61]{Al}, and Table \ref{induced} (where the last row indicates $\rho_{\xi',1}^L$ for $\xi'\in\cO$), we see that $\gamma_{\Lg^*}(\chi_{8,2})=(\widetilde{A}_2,\bar\bQ_l)$, $\gamma_{\Lg^*}(\chi_{8,4})=({A}_2,\bar\bQ_l)$, $\gamma_{\Lg^*}(\chi_{9,2})=(B_2,\bar\bQ_l)$, $\gamma_{\Lg^*}(\chi_{12,1})=(F_4(a_3),\bar\bQ_l)$ and $\gamma_{\Lg^*}(\chi_{8,1})=(B_3,\bar\bQ_l)$.

Note that the representative $\xi$ in Table \ref{tab:F4 coadjoint} for the orbit $A_2+\widetilde{A}_1$ (resp.  $(B_3)_2$, $C_3$, $(A_2)_2$, $A_1$) is a regular nilpotent element in $\Ll^*$ for a Levi $L$ of type $A_2\widetilde{A_1}$ (resp.  $B_3$, $C_3$, $A_2$, $A_1$). Thus using \eqref{b function and Springer fiber}, \eqref{Levi} and again \cite[Tables 60-61]{Al}, we see that $\gamma_{\Lg^*}^{-1}(\cO_\xi,\bar\bQ_l)=\chi_{4,3}$  (resp. $\chi_{2,1},\ \chi_{8,3},\ \chi_{1,2},\ \chi_{2,4}$). 
For the orbit $(A_2)_2$, we have used that  $\on{Ind}_{W(A_2)}^W(3)=\on{Ind}_{W(B_3)}^W([3:-]\oplus[2:1]\oplus[1:2]\oplus[-:3])$, and for the orbit $(A_1)$, we have used that $\on{Ind}_{W(A_1)}^W(2)=\on{Ind}_{W(B_3)}^W([3:-]\oplus[2:1]\oplus[1:2]\oplus[-:3]\oplus[2\,1:-]\oplus[2:1]\oplus[1:2]\oplus[-:2\,1]\oplus[1^2:1]\oplus[1:1^2])$. 

For the orbit $\widetilde{A}_2+A_1$, again the representative $\xi$ is a regular nilpotent element in $\Ll^*$ for a Levi $L$ of type $\widetilde{A}_2A_1$. Thus $\gamma_{\Lg^*}^{-1}(\cO_\xi,\bar\bQ_l)\in\{\chi_{6,1},\chi_{9,3}\}$. Note that $\langle [21:-],\on{Res}^W_{W(B_3)}\chi_{9,3}\rangle\neq 0$ while $\langle [2:1],\on{Res}^W_{W(B_3)}\chi_{9,3}\rangle= 0$. This contradicts with \eqref{restriction} as $[2:1]=(\gamma_{\Lg^*}^{B_3})^{-1}((2;1^2_1),\bar\bQ_l)$ and $[21:-]=(\gamma_{\Lg^*}^{B_3})^{-1}((2;1^2_1),\cE_{(1^2)})$. Hence $\gamma_{\Lg^*}^{-1}(\widetilde{A}_2+A_1,\bar\bQ_l)=\chi_{6,1}$.

Let  $\xi$ be in the orbit $F_4(a_1)$. Taking $\xi'$ varying in $\Ll^*$ for a Levi subgroup $L$ of type $B_3$ and using \cite[Table 60]{Al}, we see that $n_{\xi,\xi',1,1}\neq 0$ if and only if $\xi'$ satisfies $\rho_{\xi',1}^L=[3:-]$  or $\rho_{\xi',1}^L=[2:1]$. In the former case $A_L(\xi')=1$ and in the latter case we have $A_L(\xi')=S_2$ and $\rho_{\xi',(1^2)}^L=[2\,1:-]$ (see \cite{X3}). By \eqref{restriction}, $\langle\on{Res}^W_{W(B_3)}\rho_{\xi,(1^2)},\rho'\rangle_{W(B_3)}=0$ for any $\rho'\in\on{Irr}(W(B_3))-\{[3:-],[2:1],[21:-]\}$. Thus $\rho_{\xi,(1^2)}=\chi_{2,3}$. The same argument shows that for $\xi$ in the orbit $B_2$, $\langle\on{Res}^W_{W(B_3)}\rho_{\xi,(1^2)},\rho'\rangle_{W(B_3)}=0$ for any $\rho'\in\on{Irr}(W(B_3))-\{[2:1],[21:-],[1:2],[-:3],[-:21]\}$, thus
$\rho_{\xi,(1^2)}=\chi_{4,1}$. 

Consider now $\xi_6$ in the orbit $F_4(a_3)$ and $1\neq\phi\in\on{Irr}(A_G(\xi_6))$. The same argument as above shows that $\rho_{\xi_6,\phi}\neq\chi_{2,2}$. Similarly, for $\xi_8$ in the orbit $C_3(a_1)$, $\rho_{\xi_8,(1^2)}\neq\chi_{2,2}$. This forces that $\chi_{2,2}=\gamma_{\Lg^*}(\widetilde{A}_1,\cE_{(1^2)})$.

It remains to determine the images of $\chi_{9,3},\,\chi_{6,2},\,\chi_{4,4}$ and $\chi_{1,3}$ under $\gamma_{\Lg^*}$. We make use of \eqref{restriction of component gp}. Let $\xi=e_{pqr}'+e_{qrs}'+e_{pq2r}'+e_{q2r2s}$ in the orbit $F_4(a_3)$ and $\xi'=e_{pqr}'+e_{pq2r}'$. Let $P$ be the standard parabolic subgroup $P_J$, $J=\{p,q,r\}$ and $L$ the standard Levi subgroup of $P$. Then $L$ is of type $B_3$. Moreover, $\xi'\in\Ll^*$ and $\xi\in\xi'+\Ln_P^*$. It has been shown that the orbit of $\xi$ is induced from that of $\xi'$ (see \eqref{ex-ind}).  Using the description of the component group $A_G(\xi)$ in \eqref{component group s4}, one can check that  the groups defined in \eqref{groups N H} are as follows
\beqn
\text{$H=Z_P(\xi)/Z_G^0(\xi)\cong\langle\gamma_1,\gamma_3\rangle\cong S_2\times S_2$ and $N=Z_{Z_L^0(\xi')U_P}(\xi)/Z_G^0(\xi)\cong\langle\gamma_1\rangle\cong S_2$.}
\eeqn
One checks that
\beq\label{Induced reps}
\on{Ind}_H^{A_G(\xi)}\mathbf{1}=\on{Ind}_{S_2\times S_2}^{S_4}\mathbf{1}=(4)\oplus(3,1)\oplus(2,2),\ \ \on{Ind}_H^{A_G(\xi)}\tilde{\phi'}=(3,1)\oplus(2,1,1),
\eeq
where $\mathbf{1}$ is the trivial representation and $\tilde{\phi'}$ is the lift of the sign character $\phi'$ of $H/N\cong S_2$ to $H\cong S_2\times S_2$. It is shown in \cite{X3} that 
\beqn
\text{$\rho_{\xi',1}^L=[1:1^2]$ and $\rho_{\xi',\phi'}^L=[1^3:-]$. }
\eeqn
Let $\phi\in \on{Irr}A_G(\xi)$. Using \eqref{restriction formula}, \eqref{restriction of component gp} and \eqref{Induced reps}, we conclude that
\beqn
\langle\rho_{\xi,\phi}^G,\on{Ind}_{W(B_3)}^W[1;1^2]\rangle_W=\langle\rho_{\xi,\phi}^G,\chi_{12,1}\oplus\chi_{6,2}\oplus\chi_{9,3}\rangle_W=\langle\phi,(4)\oplus(3,1)\oplus(2,2)\rangle_{S_4}
\eeqn
\beqn
\langle\rho_{\xi,\phi}^G,\on{Ind}_{W(B_3)}^W[1^3;0]\rangle_W=\langle\rho_{\xi,\phi}^G,\chi_{4,4}\oplus\chi_{9,3}\oplus\chi_{1,3}\rangle_W=\langle\phi,(3,1)\oplus(2,1,1)\rangle_{S_4}
\eeqn
where in the first identity we have used \cite[Table 60]{Al}. It thus follows that
\beq\label{characters}
\rho_{\xi,(3,1)}=\chi_{9,3},\ \rho_{\xi,(2,2)}=\chi_{6,2},\ \rho_{\xi,(2,1,1)}\in\{\chi_{4,4},\chi_{1,3}\}.
\eeq

Now let $P$ be the standard parabolic subgroup $P_J$, where $J=\{p,q,s\}$, and let $L$ be the standard Levi subgroup of $P$, which is of type $A_2\widetilde{A}_1$. Then $\xi\in\Ln_P^*$ and $\xi$ is induced from the zero orbit in $\Ll^*$. Again using \eqref{component group s4} we see that in this case
\beqn
\text{$H=Z_P(\xi)/Z_G^0(\xi)\cong\langle\gamma_1,\gamma_2\rangle=N\cong S_3$.}
\eeqn
One checks that
\beqn
\on{Ind}_H^{A_G(\xi)}\mathbf{1}=\on{Ind}_{S_3}^{S_4}\mathbf{1}=(4)\oplus(3,1).
\eeqn
The same argument as above shows that
\beqn
\langle\rho_{\xi,\phi}^G,\chi_{4,4}\oplus\chi_{12,1}\oplus\chi_{9,3}\rangle_W=\langle\phi,(4)\oplus(3,1)\rangle_{S_4}.
\eeqn
In view of \eqref{characters}, it follows from the above equation that
\beqn
\rho_{\xi,(2,1,1)}=\chi_{1,3}\text{ and }\chi_{4,4}\neq \rho_{\xi,(1^4)}.
\eeqn
The last inequality forces that $\chi_{4,4}=\rho_{\xi_8,(1^2)}$, where $\xi_8$ is in the orbit $C_3(a_1)$. This completes the verification of Table \ref{tab:F4SCg*}.

\subsection{A set of Weyl group representations}In this subsection we define a subset of $\on{Irr}(W)$ following \cite{Lu1} which gives an a prior definition of the set $\{\gamma_{\Lg^*}^{-1}(\cO,\bar\bQ_l)\,|\,(\cO,\bar\bQ_l)\in\cA_{\Lg^*}\}$. 

Let $R^\vee$ be the set of coroots. For $\alpha\in R$, let $\alpha^\vee$ denote the corresponding coroot. Define
\beqn
\widetilde{\Theta}=\{\beta\in R\,|\,\beta^\vee-\alpha^\vee\notin R^\vee,\forall\ \alpha\in\Pi\};
\eeqn
\beqn
\begin{gathered}
\widetilde{\Theta}_r=\{J\subset\widetilde{\Theta}\,|\,J\text{ linearly independent, }|\sum_{\alpha\in\Pi}\mathbb{Z}\alpha^\vee/\sum_{\beta\in J}\mathbb{Z}\beta^\vee|=r^k\text{ for some }k\in\mathbb{Z}_{>0}\}.
\end{gathered}
\eeqn
For $J\in\widetilde{\Theta}_r$, let $W_J$ be the subgroup of $W$ generated by the reflections $s_\alpha,\alpha\in J$.

We define a set
$\mathcal{T}_W^{*,r}\subset\on{Irr}W$ by induction on $|W|$ as follows. If
$W=\{1\}$, $\mathcal{T}_W^{*,r}=\on{Irr}(W)$. If $W\neq\{1\}$, then
\beqn
\mathcal{T}_W^{*,r}= \mathcal{S}_W^1\cup\{\rho\in \on{Irr}(W)\,|\,  \rho=j_{W_J}^W\rho'\text{ for some
}J\in\widetilde{\Theta}_r\text{ and some }\rho'\in\mathcal{T}^{*,r}_{ W_J}\},
\eeqn 
where $\cS^1_W$ is defined as in \cite[1.3]{Lu2}, i.e. it is the set of irreducible characters of $W$ that correspond to the pairs $(\cO,\bar\bQ_l)$ in characteristic 0. 

One can verify  that the set 
\beqn
\begin{gathered}
\text{$\{\gamma_{\Lg^*}^{-1}(\cO,\bar\bQ_l)\,|\,(\cO,\bar\bQ_l)\in\cA_{\Lg^*}\}$ coincides with the set $\mathcal{T}_W^{*,3}$  (resp. $\mathcal{T}_W^{*,2}$)}\\\text{ for $G_2$ (resp. $F_4$) in characteristic 3 (resp. 2).}
\end{gathered}
\eeqn
In fact for $G_2$, we have $\widetilde{\Theta}_3=\{J\}$, where $J=\{\alpha,-(2\alpha+\beta)\}$ and $W_J=W(A_2)$, and  in characteristic 3 
\beqn
\{\gamma_{\Lg^*}^{-1}(\cO,\bar\bQ_l)\,|\,(\cO,\bar\bQ_l)\in\cA_{\Lg^*}\}=\mathcal{S}_W^1=\mathcal{T}_W^{*,3};
\eeqn for $F_4$, we have $\widetilde{\Theta}_2=\{J_1,\,J_2,\,J_3\}$, where $J_1=\{p,q,r,-p2q3r2s\}$, $W_{J_1}=W(B_3A_1)$, $J_2=\{p,r,s,-p2q3r2s\}$, $W_{J_2}=W(\widetilde{A}_3A_1)$, $J_3=\{q,r,s,-p2q3r2s\}$ and $W(J_4)=W(C_4)$, and in characteristic 2
\beqn
\{\gamma_{\Lg^*}^{-1}(\cO,\bar\bQ_l)\,|\,(\cO,\bar\bQ_l)\in\cA_{\Lg^*}\}=\mathcal{S}_W^1\cup\{\chi_{2,1},\,\chi_{1,2}\}=\mathcal{T}_W^{*,2}.
\eeqn 
Here we have used again the tables in \cite{Al} and the description of $\mathcal{T}_W^{*,2}$ in \cite{X3} when $W$ is of type $B$ or $C$.

\section{Springer correspondence for $\Lg$}\label{sec-SCg}
Let $G$ be of type $G_2$ (resp. $F_4$) defined over an algebraically closed field of characteristic $2$ or $3$ (resp.  characteristic $3$). In this section, we describe the Springer correspondence maps $\gamma_\Lg:\on{Irr}W\to\cA_\Lg$ explicitly, see Tables \ref{tab:G2 SCg char2}, \ref{tab:G2 SCg char3} and \ref{tab:F4 SCg char3}. In particular, as mentioned in the introduction, the Springer correspondence for $\Lg$ turns out to be the same as in characteristic $0$ when $G$ is of type $G_2$ and $\on{char}\tk=2$, or when $G$ is of type $F_4$ and $\on{char}\tk=3$. In view of Lemma \ref{bilinear} and \S\ref{sec-sc}, this completes the description of the Springer correspondence for $\Lg$ and $\Lg^*$ when $G$ is of type $G_2$ and $F_4$ in all characteristics, as explained in the introduction.

As in \S\ref{sec-sc}, the orbits are labelled as in the Bala-Carter classification, and $\rho_{x,\phi}$, for $x\in \cN_{\Lg}$ and $\phi\in\on{Irr}A_G(x)$, denotes the inverse image of $(\cO_x,\cE)$ under $\gamma_\Lg$, where $\cE$ is the local system on the orbit $\cO_x$ of $x$ corresponding to $\phi$.

\begin{table}[h]
\centering
\bgroup
\def\arraystretch{1.3}
\begin{tabular}{c c c c c c c}
\hline
Orbit of $x$&$\phi$&$\rho_{x,\phi}$& \vline\ \vline& Orbit of $x$ &$\phi$ &$\rho_{x,\phi}$\\
\hline
$G_2$&$(1)$&$\chi_{1,1}$& \vline\ \vline&$\widetilde{A}_1$&$(1)$&$\chi_{2,2}$\\
$G_2(a_1)$&$(3)$&$\chi_{2,1}$& \vline\ \vline&${A_1}$&$(1)$&$\chi_{1,4}$\\
$G_2(a_1)$&$(2,1)$&$\chi_{1,3}$&\vline\ \vline&
$\emptyset$ &(1)&$\chi_{1,2}$ \\
$G_2(a_1)$&$(1^3)$&$-$&\vline\ \vline&\\
\hline
\end{tabular}
\egroup
\vspace{.1in}
\caption{Springer correspondence for $\Lg$, type $G_2$ and $\on{char}\tk=2$}
\label{tab:G2 SCg char2}
\end{table}
\begin{table}[h]
\centering
\bgroup
\def\arraystretch{1.3}
\begin{tabular}{c c c c c c c}
\hline
Orbit of $x$&$\phi$&$\rho_{x,\phi}$& \vline\ \vline& Orbit of $x$ &$\phi$ &$\rho_{x,\phi}$\\
\hline
$G_2$&$(1)$&$\chi_{1,1}$& \vline\ \vline&$(\widetilde{A}_1)_2$&$(1)$&$\chi_{1,3}$\\
$G_2(a_1)$&$(2)$&$\chi_{2,1}$& \vline\ \vline&${A_1}$&$(1)$&$\chi_{1,4}$\\
$G_2(a_1)$&$(1,1)$&$-$&\vline\ \vline&
$\emptyset$ &(1)&$\chi_{1,2}$ \\
$\widetilde{A}_1$&$(1)$&$\chi_{2,2}$&\vline\ \vline&\\
\hline
\end{tabular}
\egroup
\vspace{.1in}
\caption{Springer correspondence for $\Lg$, type $G_2$ and $\on{char}\tk=3$}
\label{tab:G2 SCg char3}
\end{table}

\begin{table}[h]
\centering
\bgroup
\def\arraystretch{1.3}
\begin{tabular}{c c c c c c c c c c cccc}
\hline
Orbit of $x$ &$\phi$ &$\rho_{x,\phi}$& \vline\ \vline& Orbit of $x$ &$\phi$ &$\rho_{x,\phi}$& \vline\ \vline& Orbit of $x$ &$\phi$ &$\rho_{x,\phi}$\\
\hline
$F_4$&$(1)$&$\chi_{1,1}$& \vline\ \vline&$F_4(a_3)$&$(2^2)$& ${\chi_{6,2}}$& \vline\ \vline&$\widetilde{A_2}$&$(1)$&$\chi_{8,2}$\\
$F_4(a_1)$&$(2)$&$\chi_{4,2}$& \vline\ \vline&$F_4(a_3)$&$(2,1,1)$& ${\chi_{1,3}}$& \vline\ \vline&$A_2$&$(2)$&$\chi_{8,4}$ \\
$F_4(a_1)$&$(1^2)$ &$\chi_{2,3}$& \vline\ \vline&$F_4(a_3)$&$(1^4)$& $-$&\vline\ \vline&$A_2$&$(1^2)$&$\chi_{1,2}$ \\
$F_4(a_2)$&$(2)$&$\chi_{9,1}$&\vline\ \vline&$C_3(a_1)$&$(2)$&$\chi_{16,1}$&\vline\ \vline&$A_1+\widetilde{A_1}$&$(1)$&$\chi_{9,4}$ \\
$F_4(a_2)$&$(1^2)$&$\chi_{2,1}$& \vline\ \vline&$C_3(a_1)$&$(1^2)$ &${\chi_{4,4}}$&\vline\ \vline&$\widetilde{A_1}$&$(2)$&$\chi_{4,5}$\\
$B_3$&$(1)$&$\chi_{8,1}$&\vline\ \vline&$B_2$&$(2)$&$\chi_{9,2}$& \vline\ \vline&$\widetilde{A_1}$&$(1^2)$ &$\chi_{2,2}$\\
$C_3$&$(1)$&$\chi_{8,3}$&\vline\ \vline&$B_2$&$(1^2)$ &$\chi_{4,1}$&\vline\ \vline&$A_1$&$(1)$&$\chi_{2,4}$\\
$F_4(a_3)$&$(4)$&$\chi_{12,1}$&\vline\ \vline&$\widetilde{A_2}+A_1$&$(1)$&$\chi_{6,1}$&\vline\ \vline&$\emptyset$&$(1)$&$\chi_{1,4}$\\
$F_4(a_3)$&$(3,1)$& ${\chi_{9,3}}$& \vline\ \vline&$A_2+\widetilde{A_1}$&$(1)$&$\chi_{4,3}$&\vline\ \vline&\\
\hline
\vspace{.1in}
\end{tabular}
\egroup
\caption{Springer correspondence for $\Lg$, type $F_4$ and $\on{char}\tk=3$}
\label{tab:F4 SCg char3}
\end{table}

\subsection{Nilpotent orbits in $\Lg$}\label{ssec-nilp orbit}

The nilpotent orbits in $\Lg$ for $G_2$ have been classified by Stuhler \cite{St} and those for $F_4$ by Holt and Spaltenstein \cite{HS}. It has been shown (see \cite{HS}) that for $x\in\cN_\Lg$,
\beqn
\dim\cB_x=\frac{\dim Z_G(x)-\on{rank}G}{2},
\eeqn
where $\cB_x=\{gB\,|\,\on{Ad}(g^{-1})x\in\Lb\}$. 
We list the results here for the reader's convenience, see Tables \ref{tab:G2 g char2}, \ref{tab:G2 g char3} and \ref{tab:F4 g char3}.

Note that the nilpotent orbits in $\Lg$ for $F_4$ in $\on{char}\tk=3$ and component groups of centralizers coincide with those in characteristic $0$. In particular, there are $4$ distinguished nilpotent orbits (where $x\in\cN_\Lg$ is distinguished if $Z(G)^0$ is a maximal torus in $Z_G(x)$) and every orbit contains an element that is distinguished in a Levi subgroup. \begin{table}[h]
\centering
\bgroup
\def\arraystretch{1.2}
\begin{tabular}{c|c|c|c|c }
\hline
Orbit&representative $x$&$\dim Z_G(x)$&$A_G(x)$&$\dim\cB_x$\\
\hline
$G_2$&$e_\alpha+e_\beta$&$2$&1&0\\
$G_2(a_1)$&$e_\beta+e_{2\alpha+\beta}$&$4$&$S_3$&1\\
$\widetilde{A_1}$&$e_\alpha$&$6$&1&2\\
${A_1}$&$e_\beta$&$8$&1&3\\
$\emptyset$ &0&$14$&1&6\\
\hline
\end{tabular}
\egroup
\vspace{.1in}
\caption{Type $G_2$, nilpotent orbits, $\on{char}\tk=2$}
\label{tab:G2 g char2}
\end{table}
\begin{table}[h]
\centering
\bgroup
\def\arraystretch{1.2}
\begin{tabular}{c|c|c|c|c }
\hline
Orbit&representative $x$&$\dim Z_G(x)$&$A_G(x)$&$\dim\cB_x$\\
\hline
$G_2$&$e_\alpha+e_\beta$&$2$&1&0\\
$G_2(a_1)$&$e_\beta+e_{2\alpha+\beta}$&$4$&$S_2$&1\\
$\widetilde{A_1}$&$e_{\alpha+\beta}+e_\beta$&$6$&1&2\\
$(\widetilde{A_1})_2$&$e_\alpha$&$8$&1&3\\
${A_1}$&$e_\beta$&$8$&1&3\\
$\emptyset$ &0&$14$&1&6\\
\hline
\end{tabular}
\egroup
\vspace{.1in}
\caption{Type $G_2$, nilpotent orbits, $\on{char}\tk=3$}
\label{tab:G2 g char3}
\end{table}
\begin{table}[h]
\centering
\bgroup
\def\arraystretch{1.2}
\begin{tabular}{c|c|c|c|c}
\hline
Orbit&representative $x$&$\dim Z_G(x)$&$A_G(x)$&$\dim \cB_x$\\
\hline
$F_4$&$e_p+e_q+e_r+e_s$&$4$&1&0\\
$F_4(a_1)$&$e_p+e_{qr}+e_{q2r}+e_s$&$6$&$S_2$&1\\
$F_4(a_2)$&$e_p+e_{qr}+e_{rs}+e_{q2r2s}$&$8$&$S_2$&2\\
$B_3$&$e_p+e_{qr}+e_{q2r2s}$&${10}$&1&3\\
$C_3$&$e_s+e_{q2r}+e_{pqr}$&$10$&1&3\\
$F_4(a_3)$&$e_{pqr}+e_{qrs}+e_{pq2r}+e_{q2r2s}$&${12}$&$S_4$&4\\
$C_3(a_1)$&$e_{pqr}+e_{q2rs}+e_{q2r2s}$&$14$&$S_2$&5\\
$B_2$&$e_{pqr}+e_{q2r2s}$&$16$&$S_2$&6\\
$\widetilde{A_2}+A_1$&$e_{pqr}+e_{q2rs}+e_{p2q2r2s}$&$16$&1&6\\
$A_2+\widetilde{A_1}$&$e_{p2q2r}+e_{q2r2s}+e_{pq2rs}$&$18$&1&7\\
$\widetilde{A_2}$&$e_{pqrs}+e_{q2rs}$&$22$&1&9\\
$A_2$&$e_{p2q2r}+e_{pq2r2s}$&$22$&$S_2$&9\\
$A_1+\widetilde{A_1}$&$e_{p2q2r2s}+e_{p2q3rs}$&$24$&1&10\\
$\widetilde{A_1}$&$e_{p2q3r2s}$&$30$&$S_2$&13\\
$A_1$&$e_{2p3q4r2s}$&$36$&1&16\\
$\emptyset$&$0$&$52$&1&24\\
\hline
\end{tabular}
\egroup
\vspace{.1in}
\caption{Type $F_4$, nilpotent  orbits, $\on{char}\tk=3$}
\label{tab:F4 g char3}
\end{table}

\subsection{Springer correspondence}
We use the same strategy as in \S\ref{sec-sc} to obtain the explicit description of the maps $\gamma_\Lg$, following the methods in \cite{AL,S4,S3}. As the arguments are entirely similar, we omit the details and only explain one most complicated case.

For $F_4$, we make use of Table \ref{induced np} on induced nilpotent orbits. Here $\tilde{\cO}=\on{Ind}_{\Ll}^{\Lg}\cO$, $0$ denotes the zero orbit in the third row, and the nilpotent orbits in $\Ll$ correspond to a partition.
\begin{table}[h]
\begin{center}
\bgroup
\def\arraystretch{1.3}
\begin{tabular}{ >{$}l<{$} | >{$}l<{$} | >{$}l<{$} | >{$}l<{$} | >{$}l<{$} | >{$}l<{$} | >{$}l<{$} | >{$}l<{$} | >{$}l<{$}  }
\widetilde\cO&\widetilde{A}_2&A_2&B_2&F_4(a_3)&F_4(a_3)&F_4(a_1)&F_4(a_2)\\\hline
\text{type of }L&B_3&C_3&C_3&\widetilde{A}_2A_1&B_3&B_3&A_1\widetilde{A}_1\\\hline
\cO&0&0&2^11^4&0&3^11^4&5^11^2&0\\
\hline
\rho_{x',1}^L&[-:1^3]&[-:1^3]&[1^3:-]&[1^3]\times[1^2]&[1:1^2]&[2:1]&[2]\times[2]
\end{tabular}
\egroup
\end{center}
\vspace{.1in}
\caption{Induced nilpotent orbits, type $F_4$, $\on{char}\tk=3$}
\label{induced np}
\end{table}

The most complicated case is the orbit $F_4(a_3)$. Let $x=e_{pqr}+e_{qrs}+e_{pq2r}+e_{q2r2s}$ be an element in the orbit $F_4(a_3)$. By direct computations using \eqref{commutator relations-l} and \eqref{commutator relations-g}, one verifies that 
\beqn
A_G(x)\cong\langle\gamma_1,\gamma_2,\gamma_3\rangle\cong S_4,
\eeqn 
where
\beqn
\begin{gathered}
\gamma_1=h_p(-1)h_q(-1)h_r(-1)h_s(-1)x_p(-1)x_s(-1)x_{rs}(1),\\
\gamma_2=h_p(-1)h_q(-1)h_r(-1)h_s(-1)\,n_p\,n_s,\\
\gamma_3=h_p(-1)h_q(-1)x_p(-1)x_s(1)x_{rs}(1)x_r(-1)\,n_r\,x_r(-1)
\end{gathered}
\eeqn
with
\beqn
\gamma_1^2=\gamma_2^2=\gamma_3^2=1,\ (\gamma_1\gamma_2)^3=(\gamma_2\gamma_3)^3=1,\ (\gamma_1\gamma_3)^2=1.
\eeqn
In particular we have used the following formula 
\beqn
\begin{gathered}
\on{Ad}(u(\mathbf{t}))x=e_{pqr}+e_{qrs}+(1-2t_2)e_{pq2r}+t_2\,e_{q2rs}+(t_1-t_3)\,e_{pqrs}+(1-2t_4)e_{q2r2s}\\
+(t_1t_2-t_2t_3+t_3+t_4)e_{pq2rs}+(t_1+t_3^2-2t_4(t_1-t_3))e_{pq2r2s}
\end{gathered}
\eeqn
where $u(\mathbf{t})=x_p(t_1)\,x_r(t_2)\,x_s(t_3)\,x_{rs}(t_4)$, $t_i\in\mathbb{G}_a$.
One can then proceed the same way as in the case of orbit $F_4(a_3)$ in $\Lg^*$ in characteristic $2$.

\end{document}